\documentstyle[12pt]{article}
\def\beq{\begin{equation}}
\def\eeq{\end{equation}}
\def\bea{\begin{eqnarray}}
\def\eea{\end{eqnarray}}
\def\nn{\nonumber}
\def\ba{\begin{array}}
\def\ea{\end{array}}
\def\one{1\hskip -1mm{\rm l}}

\begin{document} 
\begin{flushright}
q-alg/9511028
\end{flushright}
\rightline{TIFR/TH/95-54}
\rightline{IMSc/95/29}
\rightline{November 1995}
\baselineskip14pt
\begin{center} 
{\bf MULTIPARAMETRIC AND COLOURED EXTENSIONS OF THE QUANTUM GROUP $GL_q(N)$ 
AND THE YANGIAN ALGEBRA $Y(gl_N)$ THROUGH A SYMMETRY TRANSFORMATION OF THE 
YANG-BAXTER EQUATION}\\ 

\bigskip

B. BASU-MALLICK\footnote{E-mail~: biru@theory.tifr.res.in} \\
{\em Theoretical Physics Group, Tata Institute of Fundamental Research, \\ 
Homi Bhabha Road, Bombay-400 005, India} \\

\medskip 

P. RAMADEVI\footnote{E-mail~: rama@imsc.ernet.in},~~
R. JAGANNATHAN\footnote{E-mail~: jagan@imsc.ernet.in} \\
{\em The Institute of Mathematical Sciences, \\

C.I.T.Campus, Tharamani, Madras-600 113, India} \\
\end{center} 

\vspace{2cm} 

\noindent 
{\small Inspired by Reshetikhin's twisting procedure to obtain 
multiparametric extensions of a Hopf algebra, a general `symmetry  
transformation' of the `particle conserving' $R$-matrix is found such that the 
resulting multiparametric $R$-matrix, with a spectral parameter as well as  
a colour parameter, is also a solution of the Yang-Baxter equation (YBE).  
The corresponding  transformation of the quantum YBE reveals a new relation 
between the associated quantized algebra and its multiparametric deformation. 
As applications of this general relation to some particular  cases, 
multiparametric and coloured extensions of the quantum group $GL_q(N)$ and the  
Yangian algebra $Y(gl_N)$ are investigated and their explicit realizations are 
also discussed.  Possible interesting physical applications of such extended 
Yangian algebras are indicated.}

\vspace{1.5cm}

\begin{center}
To appear in {\it Int. J. Mod. Phys. A}
\end{center}

\newpage 
 
\noindent 
{\bf 1. Introduction}  
\renewcommand{\theequation}{1.{\arabic{equation}}}
\setcounter{equation}{0}

\medskip

\noindent
The Yang-Baxter equation and its solutions have attracted much attention
in recent years due to their connection with diverse subjects like
exactly solvable systems, quantum groups, knot theory and conformal
field theory.$^{1-8}$  In particular, the universal 
${\cal R}$-matrix associated with a quantum algebra$^5$ plays a 
significant role in generating solutions of the Yang-Baxter equation (YBE)  
\beq
R_{12}(\lambda , \mu ) R_{13} (\lambda , \nu ) 
R_{23}(\mu , \nu )  =  
R_{23}(\mu , \nu )  R_{13} (\lambda , \nu ) 
R_{12}(\lambda , \mu )~,
\label {o1} 
\eeq
where $\lambda , \mu ,  \nu $ are spectral parameters, elements of 
$R(\lambda , \mu )$-matrix are c-numbers and the  standard notations, like 
$R_{12}(\lambda , \mu ) = R(\lambda , \mu ) \otimes \one $,  are used.  
Furthermore, by making some `twisting transformation' on the related
quantum algebra structure$^9$, one can also generate multiparametric, as 
well as coloured, solutions of YBE.$^{9-11}$  Such a twisting 
transformation on a Hopf algebra ${\cal U}$ yields a new universal matrix 
${\tilde {\cal R}}$, which is  related to the original ${\cal R}$-matrix as 
\beq 
{\tilde {\cal R}} = {\cal F}^{-1} {\cal R } {\cal F}^{-1}
\label {o2} 
\eeq
where ${\cal F} \in {\cal U} \otimes {\cal U} $ satisfies certain conditions 
and can be calculated explicitly for specific cases.$^9$ 

However, it is worth  noting  that,  in spite of the presence of many 
additional parameters,  solutions generated through the above mentioned 
`twisting transformation' often preserve some symmetries of the initial  
$R$-matrix elements. To see this through a simple example, one may consider 
the spectral parameter independent solution of YBE (\ref {o1}) 
\beq
R = q \sum_{i=1}^N e_{ii} \otimes e_{ii} + \sum_{i\neq j}
e_{ii} \otimes e_{jj} + (q-q^{-1}) \sum_{i<j }e_{ij}\otimes e_{ji}~,
\label {o3}
\eeq
which can be obtained from the fundamental representation of the universal 
${\cal R}$-matrix associated with the quantum algebra $U_q(sl(N))$. Twisting 
transformation on this $U_q(sl(N))$  algebra leads to a multiparametric 
generalization of (\ref{o3}) given by 
\beq
{\tilde R} =  q \sum_{i=1}^N e_{ii} \otimes e_{ii} + \sum_{i\neq j} 
\phi_{ij} 
e_{ii} \otimes e_{jj} + (q-q^{-1}) \sum_{i<j }e_{ij}\otimes e_{ji}~,
\label {o4}
\eeq 
where $\phi_{ij}$ are additional deformation parameters satisfying the 
condition $ \phi_{ij}\phi_{ji}=1 $.  Now, it is to be noted that both  
the solutions (\ref {o3}) and (\ref {o4}) obey a `particle conserving' 
symmetry  - i.e.,  the elements like $R^{kl}_{ij}$ and ${\tilde R }^{kl}_{ij}$ 
are nonvanishing only if the pair of `outgoing particles', $k$ and $l$, is a  
permutation of the pair of `incoming particles', $i$ and $j$.  In the present 
article we like to show that such a symmetry condition can be employed in a 
rather active way to construct new solutions (${\tilde R }(\lambda , \mu )$) 
of YBE from a given initial solution ($R (\lambda , \mu )$); this is achieved 
through a simple transformation which entails a modification of the $R$-matrix 
elements by just multiplying each of them with an appropriate factor.  This 
procedure, discussed in Sec.~2, which we call `symmetry transformation of 
YBE', is applicable even for the cases of spectral parameter dependent, 
nontriangular, initial $R$-matrix solutions.

As is well known, YBE (\ref {o1})  can be interpreted as the associativity
condition for the quantum Yang-Baxter equation (QYBE) given by
\beq 
R(\lambda ,  \mu) T_1(\lambda) T_2(\mu) = T_2(\mu) T_1(\lambda)
 R(\lambda , \mu),  
\label {o5} 
\eeq
where $ T_1(\lambda) = T(\lambda) \otimes \one~, T_2(\mu) = \one  \otimes 
T(\mu)~,$ and $T(\lambda )$ is a  matrix with operator valued elements. The 
above form of QYBE plays a central role  in the context of quantum integrable 
lattice as well as field models, Yangian algebra and quantum groups 
$^{2-5, 12-15}$.  So it is natural  to look for some transformation on both $R$ 
and $T$-matrices, which would lead to new solutions of QYBE (\ref {o5}) from a 
given initial solution.  In Sec.~2 we explore such a symmetry transformation 
of QYBE (\ref {o5}), using again the `particle conserving' restriction on the  
corresponding $R$-matrix elements in a crucial way. 

Next, we consider the applications of the above mentioned symmetry 
transformations of YBE and QYBE to some specific cases.  The symmetry  
transformation of the spectral parameter independent initial solution 
(\ref {o3}) yields the multiparametric solution (\ref {o4}) as a 
particular case.  More interesting things happen when one applies the symmetry 
transformation of QYBE to the initial $R$-matrix (\ref {o3}) and 
corresponding $T$-matrix elements.  In Sec.~3 it is shown that this 
transformation leads to a new Hopf algebra, which contains the algebra of the 
standard $GL_q(N)$ as a subalgebra.  Moreover, it turns out that all the 
generators of the multiparametric deformation of $GL(N)$ can be realized 
in terms of the generators of this new Hopf algebra in an elegant way.  

Subsequently, in Sec.~4, we focus our attention on the `particle conserving' 
rational solution of YBE given by
\beq
R(\lambda, \mu ) = 
(\lambda-\mu) \sum_{i,j=1}^N e_{ii}\otimes e_{jj} +
h \sum_{i,j=1}^N e_{ij}\otimes e_{ji}~,
\label {o6}
\eeq
where $(e_{ij})_{kl}= \delta_{ik}\delta_{jl}$.  This type of $R$-matrix is   
intimately connected with the Yangian algebra, which was recently found to 
play an important role in the analysis of some quantum spin chains with 
long-ranged interactions.$^{12-15}$   Applying the symmetry transformation
of YBE to $R(\lambda ,\mu )$-matrix (\ref {o6}) we are able to construct 
another solution of YBE, which in turn leads to a multiparametric and coloured  
extension of the standard Yangian algebra $Y(gl_N)$.  Furthermore, the symmetry 
transformation of QYBE helps us realize these deformed Yangian algebras 
through their standard counterpart. Finally, Sec.~5 contains the concluding 
remarks including hints on some possible interesting physical applications of 
the extended Yangian algebras.  

\bigskip

\noindent {\bf 2. Symmetry Transformations of YBE and QYBE}
\renewcommand{\theequation}{2.{\arabic{equation}}}
\setcounter{equation}{0}

\medskip

\noindent
In this section, our aim is to construct first a general transformation, 
which  can be performed on a `particle conserving' $R$-matrix such that the 
YBE (\ref {o1}) is still satisfied.  For this purpose, we express (\ref {o1}) 
in an elementwise form as 
\bea
   &   & \sum_{k_1,k_2,k_3} R_{i_1 i_2} ^ {k_1 k_2} (\lambda, \mu) 
R_{k_1 i_3}^{j_1 k_3} (\lambda, \nu)
R_{k_2 k_3}^{j_2 j_3} (\mu, \nu) \nn \\ 
   &   &  \qquad = \sum_{l_1,l_2,l_3} R_{l_1 l_2} ^ {j_1 j_2} (\lambda, \mu) 
R_{i_1 l_3}^{l_1 j_3} (\lambda,  \nu)
R_{i_2 i_3}^{l_2 l_3} (\mu,  \nu)~,
\label{t1}
\eea
where all indices run from $1$ to $N$.  Throughout the paper, all 
indices run from $1$ to $N$, unless otherwise stated, and the 
convention of summing over repeated indices is not used anywhere.  
Now, by exploiting the fact that for  the `particle conserving' 
case an element like $R_{ij}^{kl}$ can take non-zero value only if 
$i=k,j=l$ or $i=l,j=k$, one can write down a formal expression of 
this element as 
\beq
R_{ij}^{kl} (\lambda, \mu) = f_{ij} (\lambda, \mu) \delta_{ik} 
\delta_{jl} + g_{ij} (\lambda, \mu) \delta_{il} \delta_{jk}~,
\label{t2}
\eeq
where $f_{ij}(\lambda, \mu)$ and $g_{ij} (\lambda, \mu)$ are some yet 
undetermined functions. Substituting the above $R$-matrix element in YBE  
(\ref {t1}) and summing over its internal indices, we easily obtain the 
following set of equations: 
\renewcommand{\theequation}{2.3{\alph{equation}}}
\setcounter{equation}{0}
\bea
  &  & f_{i_1 i_2} (\lambda, \mu) f_{i_1 i_3} (\lambda, \nu) =  f_{i_1 i_3}
(\lambda, \mu) f_{i_1 i_2} (\lambda, \nu)~, \\ 
  &  & f_{i_2 i_3} (\lambda, \nu) f_{i_1 i_3} (\mu, \nu) =  f_{i_1 i_3}
(\lambda, \nu) f_{i_2 i_3} (\mu, \nu)~, \\ 
  &  & f_{i_1 i_2} (\lambda, \mu) g_{i_1 i_3} (\lambda, \nu) f_{i_2 i_1}
(\mu, \nu) \nn \\ 
  &  & \quad + g_{i_1, i_2} (\lambda, \mu) g_{i_2 i_3} (\lambda, \nu)
g_{i_1 i_2} (\mu, \nu) \nn \\ 
  &  & \qquad = f_{i_3 i_2} (\lambda, \mu) g_{i_1 i_3}
(\lambda, \nu) f_{i_2 i_3} (\mu, \nu) \nn \\ 
  &   & \qquad \quad + g_{i_2 i_3} (\lambda, \mu)
g_{i_1 i_2} (\lambda, \nu) g_{i_2 i_3} (\mu, \nu), \\ 
  &  & f_{i_1 i_2} (\lambda, \mu) g_{i_1 i_3} (\lambda, \nu) g_{i_2 i_1}
(\mu, \nu) \nn \\ 
  &   & \quad +g_{i_1 i_2} (\lambda, \mu) g_{i_2 i_3} (\lambda, \nu)
f_{i_1 i_2} (\mu, \nu) 
\nn \\ 
  &  & \qquad = g_{i_1 i_3} (\lambda, \mu) f_{i_1 i_2} 
(\lambda, \nu) g_{i_2 i_3} (\mu, \nu), \\  
  &  & g_{i_1 i_2} (\lambda, \mu) f_{i_2 i_3} (\lambda, \nu) g_{i_1 i_3} 
(\mu, \nu) \nn \\ 
  &  & \qquad = f_{i_2 i_3} (\lambda, \mu) g_{i_1 i_2} (\lambda, \nu)
g_{i_2 i_3} (\mu, \nu ) \nn \\ 
  &   & \qquad \quad + g_{i_3 i_2}  (\lambda, \mu) g_{i_1 i_3} 
(\lambda, \nu) f_{i_2 i_3} (\mu, \nu)~. 
\eea
Consequently, the $R$-matrix elements given by (\ref {t2}) can be 
treated as constituting a solution of YBE (\ref {t1}) if the functions 
$f_{ij}(\lambda, \mu)$ and $g_{ij} (\lambda, \mu)$ satisfy the set of
equations (2.3a-e). 

Next,  we consider  a  transformation on the functions 
$f_{ij}(\lambda, \mu)$ and $g_{ij} (\lambda, \mu)$ as given by  
\renewcommand{\theequation}{2.{\arabic{equation}}}
\setcounter{equation}{3}
\bea
f_{ij} (\lambda, \mu)\rightarrow {\tilde f}_{ij} (\lambda, \mu) & = & 
\phi_{ij} {u_j^{(1)}(\lambda) u_j^{(2)}
(\lambda) \over u_i^{(1)} (\mu) u_i^{(2)} (\mu) } f_{ij} (\lambda, \mu)~,
\nn \\
g_{ij} (\lambda, \mu) \rightarrow {\tilde g}_{ij} (\lambda, \mu) & = &
 {u_i^{(1)}(\lambda) u_j^{(2)}
(\lambda) \over u_i^{(1)} (\mu) u_j^{(2)} (\mu) } g_{ij} (\lambda, \mu)~,
\label {t4}
\eea
where $\phi_{ii} = 1$, $\phi_{ij} = {\phi_{ji}^{-1}}$ are $N(N-1)/2$ 
independent constant parameters and
 $u_i^{(1)} (\lambda)$, $u_i^{(2)} (\lambda)$ 
are $2N$ arbitrary (regular and nonvanishing) functions of the spectral 
parameter $\lambda $.  This transformation (\ref {t4}) leaves the whole set 
of relations (2.3a-e) invariant.  Therefore, the transformed 
${\tilde R}$-matrix given by 
\beq
{\tilde R}_{ij}^{kl} (\lambda, \mu) = {\tilde f}_{ij}(\lambda, \mu)
 \delta_{ik} \delta_{jl} + {\tilde g}_{ij} (\lambda, \mu) \delta_{il}
\delta_{jk}~, 
\label {t5}
\eeq
would satisfy YBE (\ref {t1}), provided the initial $R$-matrix (\ref {t2}) 
is a valid  solution. Due to this property of the transformation 
(\ref {t4}), or (\ref {t5}), we would naturally call it a `symmetry 
transformation' of YBE. A special case of this symmetry transformation,
corresponding to all $\phi_{ij}=1$, has been used earlier for constructing  
non-additive type $R$-matrices from the additive ones.$^{16}$   

To make the symmetry transformation more transparent, let us rewrite 
(\ref {t5}) in matrix form as  
\beq
{\tilde R} (\lambda, \mu) = F^{-1} (\lambda, \mu) R(\lambda, \mu)
{\hat F}^{-1} (\lambda, \mu)~, 
\label {t6} 
\eeq
where  
\bea
F (\lambda, \mu) & = & \sum_{i,j} \sqrt {\phi_{ji}}~{u_i^{(1)} (\mu) 
\over u_j^{(2)} (\lambda) } e_{ii} \otimes e_{jj}~, \nn \\  
{\hat F} (\lambda, \mu) & = & \sum_{i,j} \sqrt {\phi_{ji}}~{u_i^{(2)} (\mu) 
\over u_j^{(1)} (\lambda) } e_{ii} \otimes e_{jj}~, 
\label {t7}
\eea 
and the elements of $R(\lambda , \mu )$-matrix are given by (\ref {t2}).  
It is worth noting that for a particular case corresponding to
$u_i^{(1)} (\lambda)=u_i^{(2)} (\lambda)$, one gets
${\hat F} (\lambda, \mu) = F(\lambda, \mu)$, which reduces the symmetry 
transformation (\ref {t6}) to a form quite similar to the transformation 
(\ref {o2}) associated with the Reshetikhin twisting of a Hopf algebra.  If one 
further restricts to the particular case $u_i^{(1)} (\lambda) = u_i^{(2)} 
(\lambda)=1$, then, the symmetry transformation (\ref {t6}) generates the 
multiparametric solution (\ref {o4}) from the initial solution (\ref {o3}).

So far, we  have been looking at the solutions of YBE (\ref {o1}) which depends
 on  single-component spectral parameters.  We can also treat similarly the 
bicomponent spectral parameter dependent YBE defined by  
\bea
  &   &  R_{12}(\lambda , \alpha;  \mu, \beta ) 
R_{13} (\lambda ,\alpha; \nu, \gamma  ) 
R_{23}(\mu ,\beta; \nu, \gamma ) \nn \\ 
  &   & \quad \quad =   
R_{23}(\mu ,\beta; \nu, \gamma )  R_{13} (\lambda ,\alpha; \nu, \gamma ) 
R_{12}(\lambda, \alpha; \mu, \beta)~, 
\label {t8} 
\eea 
where the first components $ (\lambda , \mu , \nu )$ correspond to the usual  
spectral parameter $(\lambda )$ and the second components $ (\alpha , \beta , 
\gamma )$ refer to the colour parameter ($\alpha$).  It is possible to 
generalize the symmetry transformation (\ref {t6}), or (\ref {t5}), for 
such bicomponent spectral parameter dependent case and construct a solution 
of YBE (\ref {t8}) as 
\renewcommand{\theequation}{2.9{\alph{equation}}}
\setcounter{equation}{0}
\bea
 R (\lambda,\alpha; \mu,\beta) & = & F^{-1} ( \alpha , \beta) 
R(\lambda, \mu) {\hat F}^{-1} (\alpha , \beta)~, \\
R_{ij}^{kl} (\lambda, \alpha; \mu, \beta) & = & f_{ij} (\lambda, \alpha; \mu, 
\beta) \delta_{ik} \delta_{jl} \nn \\ 
  &  & \qquad +  g_{ij} (\lambda, \alpha; \mu, \beta) 
\delta_{il} \delta_{jk}~, 
\eea
where $f_{ij} (\lambda, \alpha; \mu, \beta) = \phi_{ij} {u_j^{(1)}(\alpha) 
u_j^{(2)} (\alpha) \over u_i^{(1)} (\beta) u_i^{(2)} (\beta) } f_{ij} 
(\lambda, \mu),$ $ g_{ij} (\lambda, \alpha; \mu, \beta) = 
{u_i^{(1)}(\alpha) u_j^{(2)} (\alpha) \over u_i^{(1)} (\beta) u_j^{(2)} 
(\beta) } g_{ij} (\lambda, \mu)$.  Note that at the  limit $\alpha = \lambda , 
$ $\beta = \mu ,$ along with the choice of notation  $ R (\lambda,\lambda ; 
\mu,\mu ) \equiv {\tilde R } (\lambda , \mu ) $, equations (2.9a,b) reduce to 
the transformations ((\ref {t6}), (\ref {t5})) associated with the 
standard YBE with a single spectral parameter.

Next we want to  address the more interesting issue concerning  the symmetry 
transformation associated with QYBE (\ref {o5}), which reads, in elementwise 
form,  
\renewcommand{\theequation}{2.{\arabic{equation}}}
\setcounter{equation}{9}
\beq
\sum_{m,n} R_{ij}^{mn} (\lambda, \mu) T_{mk} (\lambda) T_{nl} (\mu) =
\sum_{p,q} T_{jq} (\mu) T_{ip} (\lambda) R_{pq}^{kl}
(\lambda, \mu)~. 
\label {t10}
\eeq
Substituting  the  `particle conserving'  $R$-matrix (\ref {t2}) in the above 
equation and summing over its internal indices, we get a set of algebraic   
relations involving the elements of $T(\lambda )$:
\beq
U_1 + U_2 - U_3 - U_4 = 0,
\label {t11}
\eeq
where
\bea
U_1 = f_{ij} (\lambda,\mu) T_{ik} (\lambda) T_{jl} (\mu)~, \qquad 
U_2 = g_{ij} (\lambda, \mu) T_{jk}(\lambda) T_{il} (\mu)~,\nn \\
U_3 = T_{jl} (\mu) T_{ik} (\lambda) f_{kl} (\lambda, \mu)~, \qquad 
U_4 = T_{jk} (\mu) T_{il} (\lambda) g_{lk} (\lambda, \mu)~.  
\label {t12}
\eea

Subsequently, in analogy with the bicomponent spectral parameter dependent   
YBE (\ref {t8}), we write down the bicomponent spectral parameter 
dependent QYBE in elementwise form as 
\bea
  &   & \sum_{m,n}  R_{ij}^{mn} (\lambda, \alpha; \mu, \beta) 
T_{mk} (\lambda, \alpha) T_{nl} (\mu, \beta) \nn \\ 
  &   & \qquad  = \sum_{p,q} T_{jq} (\mu, \beta) T_{ip} (\lambda, \alpha) 
R_{pq}^{kl} (\lambda, \alpha; \mu, \beta)~, 
\label {t13}
\eea
and assume further that the corresponding $R(\lambda , \alpha; \mu , 
\beta)$-matrix is related to the $R(\lambda , \mu )$-matrix (\ref {t2}) 
through the transformation (2.9).  If one substitutes such 
$R(\lambda , \alpha; \mu , \beta)$-matrix in (\ref {t13}), that  
leads to another set of algebraic relations given by 
\beq
{\tilde U}_1 + {\tilde U}_2 - {\tilde U}_3 -  {\tilde U}_4 = 0~,
\label {t14}
\eeq
where the elements ${\tilde U}_r$ $(r \in [1,4])$ can be obtained from 
the elements $U_r$ in (\ref {t12}) through the  substitution:  
$T_{ij}(\lambda )$ $\longrightarrow$ $T_{ij}(\lambda , \alpha )$ , 
$T_{ij}(\mu )$ $\longrightarrow$ $T_{ij}(\mu , \beta )$, 
$f_{ij} (\lambda, \mu)$ $\longrightarrow$ 
$f_{ij} (\lambda, \alpha; \mu, \beta)$ and $g_{ij} (\lambda, \mu)$ 
$\longrightarrow$ $g_{ij} (\lambda, \alpha; \mu, \beta)$.  It is worth noting 
that the  algebraic relations (\ref {t14}) can be considered as a 
multiparametric as well as colour parameter dependent deformation  of the 
algebra (\ref {t11}).  Particular cases of these two algebras (\ref 
{t11}) and (\ref {t14}), related to specific choices of functions  
$f_{ij}(\lambda , \mu )$ and $g_{ij}(\lambda , \mu )$ satisfying (2.3), 
will be discussed in detail in the following sections.

Let us now look for a transformation, similar to the transformation (2.9), 
through which the operator valued elements $T_{ij}(\lambda , \alpha )$ 
appearing in the algebra (\ref {t14}) can be expressed through the  
elements $T_{ij}(\lambda )$ occuring in algebra (\ref {t11}) so that the 
solutions  of the bicomponent spectral parameter dependent QYBE 
(\ref {t13}) can be generated from a class of given initial solutions 
of QYBE (\ref {t10}) related to
  single-component spectral parameter~: such a  
transformation of $T$-matrix elements, coupled with the previously derived 
transformation (2.9) for $R$-matrix elements, would constitute a symmetry 
transformation of QYBE. 

It may be observed that the transformation (2.9) contains  the 
deformation parameters $\phi_{ij}$ and colour parameter dependent functions 
$u_i^{(1)}(\alpha ), u_i^{(2)}(\alpha )$, which are not present in the initial 
$R$-matrix elements (\ref {t2}).  So, to find out the analogous transformation 
for $T$-matrix elements, it is reasonable  to introduce some extra  generators, 
which do not occur in the original algebra (\ref {t11}).  Let $\tau_i$ 
$(i \in [1,N])$ be $N$ colour parameter independent generators and 
$G(\alpha)$ be a colour parameter dependent generator such that 
\renewcommand{\theequation}{2.15{\alph{equation}}}
\setcounter{equation}{0}
\bea
U_1 + U_2 - U_3  -  U_4  & = & 0~, \\
{}[\tau_i, \tau_j] =  [\tau_i, G (\alpha)] & = & [ G (\alpha), G (\beta)] = 0~, \\
\tau_i T_{jk} (\lambda) & = & c^i_{jk} (\lambda) T_{jk} (\lambda) \tau_i~, \\ 
G (\alpha) T_{jk} (\lambda) & = & d_{jk} (\lambda, \alpha) 
T_{jk}(\lambda) G (\alpha)~, 
\eea 
where $c^i_{jk} (\lambda)$ and $d_{jk} (\lambda, \alpha)$ are some yet 
undetermined $c$-number valued functions compatible with the associativity 
condition.  Then the extended algebra (2.15), defined through generators 
$ T_{ij}(\lambda ), ~\tau_i , ~G(\alpha ),$ is seen to contain the undeformed 
algebra (\ref {t11}) as a subalgebra.  Now, we propose a realization of 
deformed algebra (\ref {t14}) through the generators of algebra (2.15) as 
given by 
\renewcommand{\theequation}{2.{\arabic{equation}}}
\setcounter{equation}{15}
\beq
T_{ij} (\lambda, \alpha) = r_{ij} (\alpha) \tau_i \tau_j G^2(\alpha) T_{ij} 
(\lambda)~, 
\label {t16} 
\eeq
where $r_{ij} (\alpha)$ are also some still undetermined  $c$-number valued 
functions.  To check the validity of this ansatz,  we substitute the 
above expression  of $T_{ij} (\lambda, \alpha)$ in the deformed algebra 
(\ref {t14}).  Then, using the commutation relations (2.15b,c,d)
we can shift  the extra generators ($\tau_i, G(\alpha)$) to one side of each 
term in (\ref{t14}) and arrive at the relation  
\beq
\tau_i \tau_j \tau_k \tau_l G^2 (\alpha) G^2 (\beta) 
\left[ S_1 U_1 + S_2 U_2 - S_3 U_3 - S_4 U_4 \right ] = 0~,
\label {t17}
\eeq
where
\bea
S_1 & = & \phi_{ij}~ {u_j^{(1)} (\alpha) u_j^{(2)} (\alpha)
r_{ik} (\alpha) r_{jl} (\beta) \over u_i^{(1)} (\beta) u_i^{(2)} (\beta) 
 c^j_{ik}(\lambda) d_{ik}^2 (\lambda, \beta) c^l_{ik}(\lambda) }~, \nn \\   
S_2 & = & {u_i^{(1)} (\alpha) u_j^{(2)} (\alpha) r_{jk} (\alpha) r_{il} (\beta) 
\over u_i^{(1)} (\beta) u_j^{(2)} (\beta) 
c^i_{jk}(\lambda) d^2_{jk} (\lambda, \beta) c^l_{jk}(\lambda)}~, \nn \\ 
S_3 & = & \phi_{kl} {u_l^{(1)} (\alpha) u_l^{(2)} (\alpha) 
r_{ik} (\alpha)  r_{jl} (\beta) \over u_k^{(1)} (\beta) u_k^{(2)} (\beta)
 c^i_{jl}(\mu) d^2_{jl} (\mu, \alpha) c^k_{jl}(\mu)}~, \nn \\  
S_4 & = & {u_l^{(1)} (\alpha) u_k^{(2)} (\alpha) r_{jk} (\beta) r_{il} (\alpha)
\over u_l^{(1)} (\beta) u_k^{(2)} (\beta) c^i_{jk}(\mu ) d^2_{jk} 
(\mu, \alpha) c^l_{jk}(\mu)}~.
\label{t18}
\eea
Comparing (\ref {t17}) and (\ref {t11}), it is evident that (\ref {t17}) would 
be automatically satisfied if one sets 
\beq
S_1 = S_2 = S_3 = S_4~. 
\label {t19}
\eeq 
So, the expression (\ref {t16}) can be treated as a  realization of the  
deformed algebra (\ref {t14}), if the yet undetermined functions  
$c^i_{jk} (\lambda)$,  $d_{jk} (\lambda, \alpha)$ and $r_{ij}(\alpha )$ satisfy 
the conditions (\ref {t19}).  Now, it follows from the condition 
$S_1=S_3$ that the functions $c^i_{jk} (\lambda)$,  $d_{jk} (\lambda, \alpha)$
can in fact be chosen to be independent of the spectral parameter $\lambda$ and 
may  be given by 
\bea
c^i_{jk} (\lambda) & = & c^i_{jk} = \sqrt {\phi_{ik}\over \phi_{ij}}~, \nn \\  
d_{jk} (\lambda, \alpha) & = & d_{jk}(\alpha) = \sqrt { u_k^{(1)} (\alpha)
u_k^{(2)} (\alpha) \over u_j^{(1)} (\alpha) u_j^{(2)} (\alpha)}~.
\label {t20}
\eea
Taking the above form of $c^i_{jk},~ d_{jk} (\alpha)$ and using further the 
condition $ S_1 = S_2 = S_4 $, we obtain
\beq 
r_{ij} (\alpha) =  {1\over \sqrt {\phi_{ij}}} ~ {u_i^{(1)} (\alpha)
\over u_j^{(1)} (\alpha)}~. 
\label {t21} 
\eeq 
Then, the realization of the elements $T_{ij}(\lambda, \alpha)$ proposed in  
(\ref {t16}) takes the form  
\beq 
 T_{ij} (\lambda, \alpha) = {1\over \sqrt {\phi_{ij}}} {u_i^{(1)} (\alpha)
\over u_j^{(1)} (\alpha)} \tau_i \tau_j G^2(\alpha) T_{ij} (\lambda)~. 
\label {t22}
\eeq 
It is thus found that there indeed exists a transformation (\ref {t22}) through 
which the operator valued elements $T_{ij}(\lambda , \alpha )$ 
appearing in the  
deformed  algebra (\ref {t14}) can be expressed through the  elements 
$ T_{ij}(\lambda ) $ of the original algebra (\ref {t11}). However, to 
achieve this, we have introduced the extra generators $\tau_i,~G(\alpha )$ 
which satisfy the commutation relations
\renewcommand{\theequation}{2.23{\alph{equation}}}
\setcounter{equation}{0}
\bea
[\tau_i, \tau_j] & = & [\tau_i, G (\alpha)] = [ G (\alpha), G (\beta)] = 0~, \\ 
\tau_i T_{jk} (\lambda) & = & \sqrt {\phi_{ik}\over \phi_{ij}}~T_{jk} (\lambda) 
\tau_i~, \\ 
G (\alpha) T_{jk} (\lambda) & = & \sqrt { u_k^{(1)} (\alpha)
u_k^{(2)} (\alpha) \over u_j^{(1)} (\alpha) u_j^{(2)} (\alpha)}~T_{jk}(\lambda) 
G (\alpha)~. 
\eea 

It may be  observed further that the transformation (\ref {t22}) can   be 
written in an elegant matrix form as  
\renewcommand{\theequation}{2.{\arabic{equation}}}
\setcounter{equation}{23}
\beq
T (\lambda, \alpha) = {\cal M}(\alpha ) T(\lambda ) {\hat {\cal M}}(\alpha ) 
\label {t24}
\eeq 
where ${\cal M}(\alpha )$ and ${\hat {\cal M }}(\alpha )$ are  diagonal 
matrices  with operator valued  elements 
\bea 
{\cal M}_{ij}(\alpha ) & = & \sqrt { u_i^{(1)} (\alpha)
\over u_i^{(2)} (\alpha) }~\tau_i G(\alpha) \delta_{ij}~, \nn \\  
{\hat {\cal M}}_{ij} (\alpha ) & = & \sqrt { u_i^{(2)} (\alpha)
\over u_i^{(1)} (\alpha) }~\tau_i G(\alpha)~ \delta_{ij}~.
\label {t25}
\eea 
For the  special case corresponding to 
 $u_i^{(1)} (\alpha )=u_i^{(2)} (\alpha )$, one evidently gets 
${\cal M}(\alpha ) = {\hat {\cal M}}(\alpha )$.  This reduces the symmetry 
transformation of $T$-matrix (\ref {t24}) to a form analogous to the 
transformation of the ${\cal R}$-matrix (\ref {o2}) associated with 
Reshetikhin's twisted Hopf algebra.

Since, by using the expressions (2.9) and (\ref {t22}) simultaneously,
 one can construct the solutions 
of the bicomponent  spectral parameter dependent 
QYBE (\ref {t13}) from a given initial solution of the standard QYBE 
(\ref {t10}), it is natural to consider these transformations of $R$ and 
$T$  as constituting a symmetry transformation of QYBE.  Moreover, at the 
particular limit $\alpha = \lambda $, $\beta = \mu $, the QYBE 
(\ref {t13}) tends to its single spectral parameter counterpart (\ref {t10}). 
So this $\alpha = \lambda $, $\beta = \mu $, limit of the symmetry 
transformation ((2.9), (\ref {t22})) can  be used to obtain  a more general 
solution of the standard QYBE (\ref {t10}) itself, from a  given initial 
solution. In another special case corresponding to the choice  
$u_i^{(1)}(\alpha ) = u_i^{(2)}(\alpha ) = 1$ for all $i$, the transformations 
(2.9) and  (\ref {t22}) become independent of colour parameter and can also be 
used to  construct new solutions of the standard QYBE. This type of symmetry  
transformation will be used shortly, in the next section, for establishing a 
link between the single-parametric and multiparametric deformations of $GL(N)$. 

\bigskip

\noindent 
{\bf 3. A New Realization of Multiparametric Deformation of $GL(N)$} 
\renewcommand{\theequation}{3.{\arabic{equation}}}
\setcounter{equation}{0}

\medskip

\noindent
As is well known, the algebra of $GL_q(N)$ is generated by $N^2$ 
elements $T_{ij}$ satisfying  the commutation relations$^{17}$ 
\bea
T_{ij} T_{ik} & = & q^{-1} T_{ik} T_{ij}~, \qquad  
T_{ik} T_{lk} = q^{-1} T_{lk} T_{ik}~, \nn \\ 
T_{ik} T_{lj} & = & T_{lj} T_{ik}~, \qquad 
[T_{ij}, T_{lk}] = (q^{-1} - q) T_{lj} T_{ik}~, 
\label {th1}
\eea 
where $i<l,~j<k$. This algebra is obtained by substituting the `particle 
conserving' $R$-matrix (\ref{o3}) in the spectral parameterless limit of QYBE 
(\ref{o5}).  One can construct a multiparametric generalization of 
$GL_q(N)$,$^{18}$ by substituting the $R$-matrix (\ref{o4}) in QYBE 
(\ref{o5})~: with $i<l,~j<k$,  
\bea 
{\tilde T}_{ij} {\tilde T}_{ik} & = & q^{-1} 
\phi_{jk} {\tilde T}_{ik} {\tilde T}_{ij}~, \qquad  
{\tilde T}_{ik} {\tilde T}_{lk} = q^{-1} \phi_{li} {\tilde T}_{lk} 
{\tilde T}_{ik}~, \nn \\
\phi_{il} {\tilde T}_{ik} {\tilde T}_{lj} & = & \phi_{kj} {\tilde T}_{lj} 
{\tilde T}_{ik}~, \quad 
\phi_{il} {\tilde T}_{ij} {\tilde T}_{lk} - \phi_{jk} {\tilde T}_{lk} 
{\tilde T}_{ij} =  (q^{-1} -q) {\tilde T}_{lj} {\tilde T}_{ik}~,  
\label {th2}
\eea 
where $\phi_{ij}$ are $ {N(N-1)/2}$ additional deformation parameters.
Evidently, at the limit $\phi_{ij} = 1$ for all $i,j$, the multiparametric 
deformed algebra (\ref {th2}) reduces to its single parameter dependent 
version (\ref {th1}).

Various aspects of the algebra of $GL_q(N)$ (\ref{th1}) have been studied quite 
extensively in the  literature. In particular, we are interested in the 
realization of the generating elements $T_{ij}$ of $GL_q(N)$, for $|q| = 1$, 
through recasting it in the Heisenberg-Weyl form  
($A_i B_j = \omega_{ij} B_j A_i~,~|\omega_{ij}| = 1)$ and subsequently using 
mutually commuting pairs of canonically conjugate quantum mechanical operators 
(boson realization), or, finite dimensional matrices.$^{19-21}$  The 
multiparametric extensions  of $GL_q(N)$ are, in general, more difficult to 
analyze and the corresponding  realizations have  been investigated only for 
some special cases like $N =2$.$^{22,23}$ So, if one can 
express the generators of the  multiparametric algebra (\ref{th2}) in terms of 
the generators of another  algebra structurally similar to the algebra of 
$GL_q(N)$ (\ref{th1}), then,  such expression should be much helpful for 
constructing the realizations of the  multiparametric deformation of $GL(N)$.  
To this end, we use the symmetry transformation (\ref{t22}).  

As has been already noted in Sec. II, the `particle conserving'  $R$-matrices 
(\ref{o3}) and  (\ref{o4}) are related  through a particular limit  of  
symmetry transformation (2.9)  corresponding to the choice 
$u_i^{(1)}(\alpha ) = u_i^{(2)}(\alpha ) = 1,$ for all $i$. Consequently, the 
pair of algebras (\ref{th1}) and (\ref{th2}) can be considered as a special
case of the pair (\ref{t11}) and (\ref{t14}), which are connected 
through the general symmetry transformation (\ref{t22}).  More explicitly, in   
the present context, the algebra 
\renewcommand{\theequation}{3.3{\alph{equation}}}
\setcounter{equation}{0}
\bea
[\tau_i,  \tau_j] & = & 0~, \qquad  
\tau_m T_{kl} =  {\sqrt {\phi_{ml} \over  \phi_{mk}}} 
~ T_{kl} \tau_m~, \\ 
T_{ij} T_{ik} & = & q^{-1} T_{ik} T_{ij}~, \qquad  T_{ik} T_{lk} = q^{-1} 
T_{lk} T_{ik}~, \\ 
T_{ik} T_{lj} & = & T_{lj} T_{ik}~, \qquad [T_{ij}, T_{lk}] =  (q^{-1} - q) 
T_{lj} T_{ik}~, 
\eea 
where  $\tau_i$ are $N$ extra generators, and $i<l,~j<k$, contains $GL_q(N)$  
as a subalgebra, and is a particular  case of 
the algebra (2.15) corresponding to the choice  $u_i^{(1)}(\alpha ) = 
u_i^{(2)}(\alpha ) = 1$ and $G(\alpha ) = \one$.  Consequently, from 
(\ref{t22}), or (\ref{t24}), we obtain the realization 
\renewcommand{\theequation}{3.{\arabic{equation}}}
\setcounter{equation}{3}
\beq 
{\tilde T}_{ij} = {1 \over {\sqrt {\phi_{ij}}}}~\tau_i \tau_j T_{ij}~,  
\quad {\rm or}~~~{\tilde T} = {\cal M}  T {\cal M}~, 
\label{th4}
\eeq
where ${\cal M}$ is the diagonal $N\times N$ matrix with operator valued 
elements given by ${\cal M}_{ij} = \tau_i \delta_{ij} $.  By using the algebra
(3.3), one can also  verify directly the validity of above realization 
(\ref{th4}).  Thus, we find that it is indeed possible
 to express the generators 
of the multiparametric  deformation of $GL(N)$ (\ref{th2}) through the 
generators of $GL_q(N)$ (\ref{th1}) and $N$ additional generators $\tau_i$ 
satisfying the commutation relations (3.3a). 
 
Though the above observations are true in general,  sometimes one might be able 
to obtain a realization of the multiparametric  deformation of $GL(N)$, using 
less than $N$ additional generators $\tau_i$, besides the generators of 
$GL_q(N)$, provided the corresponding parameters $\phi_{ij}$ satisfy certain 
constraints.   To illustrate this point, 
let us assume that all $\phi_{ij}$  can 
be written in the form: $\phi_{ij} = {r_i \over r_j}$, where $r_i$s are $N$ 
independent parameters.  Then, the structure constant $\sqrt { \phi_{ml} \over 
\phi_{mk}}$, occuring in the commutation relation (3.3a), would be independent 
of the index $m$. Consequently,  all generators $\tau_i$ occuring in the 
algebra (3.3), as well as in the realization (\ref{th4}), can be replaced 
effectively by a single generator $\tau $ satisfying the commutation relations 
$\tau T_{kl} = \sqrt { r_l \over r_k }~T_{kl}~\tau $.  Thus, in this case we  
can have a realization of the multiparametric deformed $GL(N)$ algebra (\ref 
{th2})  by augmenting the $GL_q(N)$ algebra (\ref {th1}) with just one  
additional generator $\tau $.  A realization of this type has been considered 
earlier$^{24}$ for the special case $N=2$.

Next, one may ask the interesting question whether the extended algebra (3.3)  
itself corresponds to a quantum group. To answer  this question, we take  a 
$2N \times 2N$ dimensional $T'$-matrix and  a $ 4N^2 \times 4N^2 $ dimensional  
$R'$-matrix, whose nonvanishing elements are given by 
\renewcommand{\theequation}{3.5{\alph{equation}}}
\setcounter{equation}{0}
\bea
T^{\prime} _{ij}  & = & T_{ij}~,\qquad 
T^{\prime}_{i^{\prime}, j^{\prime}} = \delta_{ij} \tau_i~, \\ 
\left (R^{\prime}\right )_{ij}^{ij} & 
= & \left (R^{\prime} \right )_{i^{\prime}, 
j^{\prime}}^{i^{\prime}, j^{\prime}} = 1 + (q-1)\delta_{ij}~, \nn \\  
\left (R^{\prime}\right )_{i^{\prime},j}^{i^{\prime},i} & = & 
\left (R^{\prime} \right )_{i, j^{\prime}}
^{i, j^{\prime}}=\sqrt {\phi_{ij}}~, \qquad 
\left ( R^{\prime} \right )_{kl}^{lk} = q-q^{-1}~, 
\eea
where  $i'=i+N ,~ j'=j+N$, $1 \leq i,j  \leq N$  and $1 \leq  k <l  \leq 2N$.  
Then, it is found that the  algebra (3.3) is generated by substituting these 
$R'$ and $T'$ matrices in the spectral-parameterless limit of QYBE (\ref{o5}).  
So, the extended algebra  (3.3) indeed defines  a quantum group, for which 
the coproduct and other Hopf algebra operations are readily obtained.  From 
the  block diagonal structure of $T'$-matrix (3.5a) it is evident that this new 
quantum group is a deformation of the group of $2N \times 2N$ matrices having 
only $\{ M_{ij}, M_{N+k,N+k} | i,j,k = 1,2, \ldots, N \}$ as nonzero elements, 
or in other words it is a deformation of the subgroup of $GL(2N)$ isomorphic 
to $GL(N) \otimes \underbrace {GL(1) \otimes \cdots \otimes GL(1)}_N $.  

One can also easily construct the noncommutative planes associated with  this 
non-semisimple quantum group.  For this purpose,  let us consider a set of  
$2N$ coordinates $x_i,y_i$ ($i \in [1,N]$) undergoing the transformation  
\renewcommand{\theequation}{3.{\arabic{equation}}}
\setcounter{equation}{5}
\beq
x_i^{\prime} = \sum_{j=1}^N T_{ij} x_j ~, \qquad 
y_i^{\prime} = \tau_i y_i.
\label {th6}
\eeq
Using the algebra (3.3), it can be readily checked  that the following two sets 
of bilinear relations remain invariant under the transformation (\ref{th6}):
\renewcommand{\theequation}{3.7{\alph{equation}}}
\setcounter{equation}{0}
\bea
x_i x_j & = & q^{-1} x_j x_i~, \quad 
x_k y_l = \sqrt {\phi_{lk}}~y_l x_k~, \quad 
[y_k,y_l] = 0~, \\  
x_i^2 & = & 0~, \quad 
x_i x_j = -q x_j x_i~, \quad 
x_k y_l = \sqrt {\phi_{lk}}~y_l x_k~, \quad 
[y_k, y_l] = 0~,  
\eea 
where $i<j$. Therefore,  the relations (3.7a) and 
(3.7b)  represent the two noncommutative planes corresponding  to the  quantum 
group (3.3).  Let us now define $X_i = x_i y_i~,~~i = 1,2, \ldots , N$.
Then, it is found that $X_i$ obey the commutation relations  
\renewcommand{\theequation}{3.8{\alph{equation}}}
\setcounter{equation}{0}
\beq 
X_i X_j =   q^{-1} \phi_{ji}~ X_j X_i~, 
\eeq 
or 
\beq  
X_i^2 = 0~, \qquad  X_i X_j = -q \phi_{ji}  X_j X_i~, 
\eeq 
depending on whether $(x_i,y_i)$ satisfy (3.7a) or (3.7b) respectively.  It  
is interesting to observe that (3.8a) and (3.8b) coincide with the commutation 
relations of the Manin $q$-plane and its exterior plane associated with the 
multiparametric deformation of $GL(N)$$^{18}$. 

Let us note that the extended algebra (3.3) can be recast in the  
Heisenberg-Weyl form in a very simple way. Suppose ${\cal T}_{ij}$ are the 
basis elements, through which the  $GL_q(N)$ algebra  (\ref{th1}) can be 
expressed in the  Heisenberg-Weyl form for unimodular values of the parameter 
$q$; explicit construction of such ${\cal T}_{ij}$ is known.$^{19,20}$  Then, 
it can be verified directly that the extended algebra (3.3) will also take the 
Heisenberg-Weyl form if ${\cal T}_{ij}$, $\tau_k$ are chosen as the basis 
elements.  Hence,  the multiparametric deformation of $GL(N)$ can also be 
realized in terms of mutually commuting pairs of canonically conjugate quantum 
mechanical operators and finite dimensional matrices using the known 
methods.$^{21,22}$   

\bigskip

\noindent{\bf 4. Multiparametric and Coloured Extensions of $Y(gl_N)$}  
\renewcommand{\theequation}{4.{\arabic{equation}}}
\setcounter{equation}{0}

\medskip

\noindent
Here our aim is to study various types of deformations of the Yangian algebra 
$Y(gl_N)$$^{5, 12-15}$ and their interrelations, using the symmetry 
transformations of YBE and QYBE.  The standard $Y(gl_N)$, with the defining 
relations  
\beq
(\lambda - \mu ) \left[ T_{ij}(\lambda) , T_{kl}(\mu) \right] = h 
\left\{ T_{kj}(\mu ) T_{il}(\lambda ) 
- T_{kj}(\lambda ) T_{il}(\mu ) \right\}~, 
\label {f1} 
\eeq 
results from the substitution of the rational $R(\lambda )$-matrix (\ref {o6}) 
in QYBE (\ref {t10}).  If one assumes the usual analyticity property of  
$T(\lambda )$ and the asymptotic condition $T(\lambda ) \longrightarrow 1$ at 
$\lambda \longrightarrow \infty$, then the  operator valued elements 
$T_{ij}(\lambda )$ can be expanded in powers of $\lambda $ as 
\beq 
T_{ij}(\lambda) = \delta_{ij} + h \sum _{n=0}^{\infty}{t_n^{ij} \over 
\lambda^{n+1}}~. 
\label {f2} 
\eeq  
Substituting  the above  expansion in (\ref {f1}) and comparing the 
coefficients of equal powers in spectral parameters on both sides, one can 
express the $Y(gl_N)$ algebra through the modes $t_n^{ij}$ as   
\renewcommand{\theequation}{4.3{\alph{equation}}}
\setcounter{equation}{0}
\bea
\left[ t_0^{ij} , t_n^{kl} \right] & = & \delta_{il} t_n^{kj} - \delta_{kj} 
t_n^{il}~, \\
\left[ t_{n+1}^{ij} , t_m^{kl} \right] - \left[ t_{n}^{ij} , t_{m+1}^{kl} 
\right] & = & h (t_m^{kj} t_n^{il} - t_n^{kj} t_m^{il})~. 
\eea 
With the help of induction procedure, it is seen that the algebra (4.3a,b) can 
also be presented equivalently as a single relation  
\bea
\left[ t_n^{ij} , t_m^{kl} \right] & = & \delta_{il} t_{n+m}^{kj} 
- \delta_{kj} t_{n+m}^{il} \nn \\ 
  &   & \qquad +  h \sum_{p=0}^{n-1} \left( t_{m+p}^{kj}   
t_{n-1-p}^{il} - t_{n-1-p}^{kj} t_{m+p}^{il} \right)~. 
\eea 
It may be noted that at the limit $h \longrightarrow 0$, (4.3c) reduces to 
a subalgebra of $gl(N)$ Kac-Moody algebra containing its  non-negative modes. 
Consequently, this Yangian algebra might be considered as some nonlinear 
deformation  of the $gl(N)$ Kac-Moody algebra through the parameter $h$. Casimir 
operators for $Y(gl_N)$  may be   obtained  by  constructing  the   
corresponding quantum determinant. Further, it turns out that all the higher 
level generators of $Y(gl_N)$ can be realized consistently through only the 
$0$th and $1$st level generators of $Y(sl_N)$, provided a few Serre 
relations are satisfied.$^{5,15}$    

Now, for constructing a multiparametric and coloured extension of $Y(gl_N)$, 
we apply the symmetry transformation of YBE (2.9) to the particle 
conserving rational $R$-matrix (\ref{o6}).  This leads to a solution of the 
bicomponent spectral parameter dependent YBE (\ref{t8}), given by
\renewcommand{\theequation}{4.{\arabic{equation}}}
\setcounter{equation}{3}
\bea 
R(\lambda, \alpha; \mu, \beta) & = & (\lambda - \mu) \sum_{i,j=1}^N u_{ij} 
(\alpha, \beta) e_{ii} \otimes e_{jj} \nn \\ 
  &  &  \qquad + h \sum_{i,j=1}^{N} v_{ij} (\alpha, \beta) e_{ij} \otimes 
e_{ji}~, 
\label {f4}
\eea 
where $u_{ij} (\alpha, \beta) = \phi_{ij} { u_j^{(1)} (\alpha)
u_j^{(2)} (\alpha) \over u_i^{(1)} (\beta) u_i^{(2)} (\beta)}$
and  $v_{ij} (\alpha, \beta) = { u_i^{(1)} (\alpha)
u_j^{(2)} (\alpha) \over u_i^{(1)} (\beta) u_j^{(2)} (\beta)}$.
Substituting this $R(\lambda , \alpha ; \mu , \beta )$-matrix in QYBE 
(\ref{t13}), we get the desired generalization of $Y(gl_N)$ (\ref {f1})~: 
\bea 
  &   &  (\lambda - \mu)\left\{ u_{ik} (\alpha, \beta) T_{ij} (\lambda, \alpha) 
T_{kl} (\mu, \beta) \right. \nn \\ 
  &   & \quad \left. - u_{jl} (\alpha, \beta) T_{kl}(\mu, \beta) 
T_{ij}(\lambda, \alpha) \right\}  \nn \\
  &   & \qquad = h \left\{ v_{lj} (\alpha, \beta) T_{kj} (\mu, \beta ) 
 T_{il}(\lambda, \alpha) \right. \nn \\ 
  &   &  \qquad \quad \left. - v_{ik} (\alpha, \beta) T_{kj} (\lambda, \alpha) 
T_{il} (\mu, \beta) \right\}~. 
\label {f5} 
\eea  
It is evident that at the particular limit $u_i^{(1)}(\alpha ) = 
u_i^{(2)}(\alpha ) = 1,$ $\phi_{ij}=1$, the above algebra (4.5) tends to its 
standard counterpart (\ref{f1}). 
  However, at another limit corresponding to the 
choice  $u_i^{(1)}(\alpha )= u_i^{(2)}(\alpha ) =1,$ but $\phi_{ij} \neq 1$,
the algebra (\ref{f5}) becomes colourless and reduces to
\bea
  &   & (\lambda - \mu) \left\{ {\phi_{ik} T_{ij} (\lambda) T_{kl} (\mu) - 
\phi_{jl} T_{kl} (\mu) T_{ij} (\lambda)}~\right\} \nn \\ 
  &   & \qquad = h \left\{ {T_{kj} (\mu) T_{il} (\lambda) - T_{kj} (\lambda) 
T_{il} (\mu)}~\right\}~. \label {f6}
\eea 
This type of multiparametric extension of the  Yangian algebra $Y(gl_N)$ 
(\ref{f1}) and related Hopf algebra properties have been explored 
earlier.$^{25}$

Next, we like to investigate whether,  in analogy with (4.3a,b,c), the 
coloured and multiparametric Yangian algebra (\ref {f5}) can be written in 
terms of the  modes of its generators.   To this end, we expand the element
$T_{ij}(\lambda , \alpha )$ in powers of the spectral parameter $\lambda $ as
\beq 
T_{ij} (\lambda, \alpha) = \tau_{ii}(\alpha) \delta_{ij} + h \sum_{n=0}^\infty 
{t_n^{ij} (\alpha) \over \lambda^{n+1}} ~, 
\label {f7}
\eeq  
where the mode operators depend on the continuous colour parameter $\alpha $.
Substituting  this mode expansion in (\ref {f5}) and equating the terms on 
both sides with the same powers of $\lambda$ and $\mu$,  we arrive at the 
following relations:
\renewcommand{\theequation}{4.8{\alph{equation}}}
\setcounter{equation}{0}
\bea
  &   & \qquad  [\tau_{ii}(\alpha), \tau_{jj}(\beta) ] =  0~, \\  
  &   & \qquad  \tau_{jj}(\alpha)  t_m^{il}(\beta) = {u_{jl}
(\alpha, \beta) 
\over u_{ji}(\alpha, \beta)}  t_m^{il} (\beta) \tau_{jj}(\alpha)~, \\  
  &   & u_{ik} (\alpha, \beta) t_0^{ij} (\alpha) t_m^{kl} (\beta) - 
u_{jl} (\alpha, \beta) t_m^{kl} (\beta) t_0^{ij}(\alpha) \nn \\ 
  &   & \qquad ~~= \delta_{il} v_{lj}(\alpha, \beta) t_m^{kj} (\beta) \tau_{ii} 
(\alpha) - \delta_{jk} v_{ik} (\alpha, \beta) \tau_{kk} (\alpha)  
t_m^{il} (\beta)~, \\ 
  &   & \{u_{ik} (\alpha, \beta) t_{n+1}^{ij}(\alpha) t_m^{kl}(\beta) 
- u_{jl} (\alpha, \beta) t_m^{kl}(\beta) t_{n+1}^{ik} (\alpha)\} \nn \\ 
  &   & \qquad -\{ u_{ik} (\alpha, \beta)~ t_n^{ij}(\alpha) t_{m+1}^{kl} 
- u_{jl}(\alpha, \beta)~ t_{m+1}^{kl}(\beta) t_n^{ij}(\alpha)\} \nn \\ 
  &   & \qquad \quad = h \{v_{lj} (\alpha, \beta)~ t_m^{kj}(\beta) t_n^{il} 
(\alpha) - v_{ik} (\alpha, \beta) t_n^{kj} (\alpha)  t_m^{il} (\beta) \}~. 
\eea 
Using the induction procedure, the last two equations (4.8c,d) can be combined 
into the single relation 
\bea 
  &   & \{ u_{ik} (\alpha, \beta) t_n^{ij}(\alpha) t_m^{kl}(\beta) 
- u_{jl} (\alpha, \beta ) t_m^{kl}(\beta) t_n^{ij}(\alpha)\} \nn \\ 
  &   & \qquad  = \{  \delta_{il} v_{lj} (\alpha, \beta) t_{m+n}^{kj}(\beta) 
\tau_{ii} (\alpha)  - \delta_{jk} v_{ik} (\alpha, \beta) \tau_{kk} (\alpha) 
t_{m+n}^{il} (\beta)\} \nn \\ 
  &   & \qquad \qquad +h \sum_{p=0}^{n-1} \left\{ v_{lj} (\alpha, \beta) 
t_{m+p}^{kj} (\beta)  t_{n-1-p}^{il} (\alpha) \right. \nn \\ 
  &   & \qquad \qquad \qquad \qquad \qquad  \left. - v_{ik} (\alpha, 
\beta) t_{n-1-p}^{kj} (\alpha)  t_{m+p}^{il} (\beta) \right\}~. 
\eea
Thus the mode operators associated with the coloured and multiparametric   
Yangian algebra (\ref {f5}) satisfy the relations (4.8).  It is 
obvious from (4.8a) and (4.8b) that in the limit $u^{(1)}_i(\alpha ) = 
u^{(2)}_i(\alpha ) = 1$, $\phi_{ij} = 1$, $\tau_{jj}(\alpha )$ can be taken to 
be the unit operator for any $j$ so that in this limit the relations (4.8c), 
(4.8d) and (4.8e) tend to their standard counterparts (4.3a), (4.3b) and 
(4.3c) respectively.    

It should  be of physical relevance to construct the representations of the  
coloured and multiparametric Yangian algebra since the representation theory 
of the standard Yangian algebra$^{26-28}$ has been found to play an important 
role in the analysis of the degeneracies of wavefunctions of some quantum 
spin chains with long-ranged interactions$^{12-15}$.  One may also ask whether 
the new mode algebra (4.8) can be realized through some $0$th and $1$st level 
generators satisfying the Serre relations.   To study these matters, we first 
like to see the connection between the deformed Yangian algebra (4.8) and its 
standard counterpart (4.3) through the symmetry transformation of QYBE 
(\ref {t22}).   To this end, let us extend the standard $Y(gl_N)$ algebra 
(4.3) by augmenting it with extra generators, $N$ colour parameter independent 
$\tau_i$ and a colour parameter dependent $G(\alpha )$, such that 
\renewcommand{\theequation}{4.9{\alph{equation}}}
\setcounter{equation}{0}
\bea
[ \tau_i, \tau_j ] & 
= & [\tau_i , G(\alpha)] = [G(\alpha) , G(\beta)] = 0~, \\  
\tau_it_n^{jk} & = & \sqrt { \phi_{ik}  \over \phi_{ij}}~t_n^{jk} \tau_i~, \\  
G(\alpha) t_n^{jk} & = & \sqrt {u_k^{(1)} (\alpha) u_k^{(2)} (\alpha) \over
u_j^{(1)} (\alpha) u_j^{(2)} (\alpha) }~t_n^{jk} G(\alpha)~, \\
\left[ t_n^{ij} , t_m^{kl} \right ] & = & \delta_{il} t_{n+m}^{kj} - 
\delta_{kj} t_{n+m}^{il} \nn \\  
  &   & \qquad + h \sum_{p=0}^{n-1} \left ( t_{m+p}^{kj}  
t_{n-1-p}^{il} - t_{n-1-p}^{kj} t_{m+p}^{il} \right )~. 
\eea 
Note that the above algebra is a special case of (2.15) and the relations 
(4.9b) and (4.9c) are obtained from (2.23b) and (2.23c), respectively, 
by inserting in them the mode expansion (\ref{f2}). Subsequently, using the 
general symmetry transformation (\ref {t22}) and the mode expansions 
(\ref {f2}) and (\ref {f7}), we find that the coloured and multiparametric 
Yangian algebra (4.8) can be realized through the generators of the extended 
algebra (4.9) as 
\renewcommand{\theequation}{4.{\arabic{equation}}}
\setcounter{equation}{9}
\beq 
\tau_{ii}(\alpha) = \left (\tau _iG(\alpha)\right )^2~, \qquad  
t_n^{ij}(\alpha) = {1\over \sqrt {\phi_{ij}}} {u_i^{(1)} (\alpha) \over 
u_j^{(1)} (\alpha)} \tau_i \tau_j G^2(\alpha) t_n^{ij}~.   
\label {f10} 
\eeq  
One can also verify the validity of this realization directly by using the 
extended algebra (4.9).

The above realization of the deformed Yangian algebra is seen to be quite 
similar to the previously constructed realization (\ref{th4}) of the 
multiparametric extension of $GL_q(N)$.  However, the extra generators 
$\tau_i$ appearing in the algebra (3.3) seem to be independent of the
original $GL_q(N)$ generators $T_{ij}$ and that is why it is difficult to 
give a realization of the algebra of multiparametric deformed $GL(N)$ 
(\ref {th2}) in terms of only the $GL_q(N)$ generators.  On the other hand, it 
turns out interestingly that one can express the extra generators $\tau_i,~
G(\alpha ) $  
appearing in the extended $Y(gl_N)$ algebra (4.9) through the original 
$Y(gl_N)$ generators, satisfying the algebra (4.3), as 
\beq  
\tau_j  = \exp \left ( i \sum_{l=1}^N A_{jl} t_0^{ll} \right )~, \qquad  
G(\alpha) = \exp \left ( i\sum_{l=1}^N B_l (\alpha) t_0^{ll} \right )~,
\label {f11} 
\eeq 
where $\exp (i A_{jl}) = \sqrt{\phi_{jl}}$ and $\exp \left( iB_l(\alpha)\right)
 = \sqrt {u_l^{(1)} (\alpha) u_l^{(2)} (\alpha) }$.  Now, substituting the 
expressions for $\tau$\,s and $G(\alpha )$ given by (\ref {f11}) in 
(\ref {f10}), we have a realization of the coloured and multiparametric Yangian 
algebra (4.8) entirely in terms of the standard Yangian generators: 
\bea 
\tau_{jj}(\alpha) & = & \exp \left (2i \sum_l [A_{jl}+B_l(\alpha)]  
t_0^{ll}\right )~, \nn \\
 t_n^{jk}(\alpha) & = & {1\over \sqrt {\phi_{jk}}} {u_j^{(1)} (\alpha) \over 
u_k^{(1)} (\alpha)} \exp \left (i \sum_{l=1}^N [A_{jl}+A_{kl}+ 2 B_l(\alpha)]
t_0^{ll}\right ) t_n^{jk}~.  
\label {f12}
\eea 
However, it should be noted that, though expressible entirely in terms of the 
generators of its standard counterpart, a multiparametric generalization of the 
Yangian algebra is endowed with a new coalgebra structure as can be seen by  
following the treatment of the Hopf algebraic structure of the algebra 
(\ref{f6}) presented in Ref. 25.  

As has been already mentioned, all the higher  level generators of the standard 
$Y(gl_N)$ (4.3) can be expressed in terms of the level-0 and level-1 generators 
of $Y(sl_N)$ provided they satisfy a few Serre relations.  Then, using the 
realization (\ref {f12}),  one can also express all higher level generators of 
the coloured and multiparametric Yangian algebra (4.8) through the level-0 
and level-1 generators of $Y(sl_N)$.  Moreover, it is also possible to  
construct the representations of the deformed Yangian algebra by 
using the realization (\ref{f12}) and the known representations of the 
standard Yangian algebra.  It is thus found that, in spite of the apparently 
complicated nature of the coloured and multiparametric Yangian algebra (4.8), 
one can easily understand it through its realization in terms of the standard 
Yangian generators. 

\bigskip

\noindent{\bf 5. Concluding Remarks}

\medskip

\noindent
The YBE and its solutions, $R$-matrices, are well known to play the central 
role in several problems associated with exactly solvable models, quantum 
algebras, knot theory, conformal field theory, etc.  So, it is useful to 
search for new solutions of the YBE.  In this article, inspired by 
Reshetikhin's twisting procedure for obtaining the multiparametric extensions 
of a Hopf algebra, we have demonstrated the operation of a `symmetry 
transformation' which, applied to a known `particle conserving' $R$-matrix, 
leads to a new spectral parameter dependent solution of the YBE depending, in 
general, on multiple $q$-parameters as well as a colour parameter.  The YBE 
follows from the associativity condition of QYBE.  Hence, naturally, we find 
that a certain set of transformations on both the $R$-matrix and the 
$T$-matrix, constituting a symmetry transformation of QYBE, leads to a 
multiparametric and coloured extension of the associated quantum group 
structure; here again, the $R$-matrix is considered to have the `particle 
conserving' symmetry.  

As an application of the symmetry transformations of YBE and QYBE, first we 
have analysed the construction of the multiparametric extension of $GL_q(N)$.  
It is found that the resulting quantum group admits a realization through the 
generators of $GL_q(N)$ and some mutually commuting extra generators.  We hope 
that such a simple realization of the multiparametric extension of $GL_q(N)$ 
would be helpful in building up the corresponding representation theory.  Next, 
by applying the symmetry transformation to a class of rational solutions of 
YBE, w have constructed a multiparametric and coloured extension of the Yangian 
algebra $Y(gl_N)$.  It is found that this extended Yangian algebra can also be 
realized completely in terms of the standard Yangian generators, indicating 
that the extra deformation occurs essentially in the coalgebra sector of the 
Hopf algebra.  

Finally, we like to hint at some possible interesting physical applications of 
the extended Yangian algebras which have been discussed here.  As is well 
known,$^{12}$ the conserved quantities of $su(N)$-invariant spin 
Calogero-Sutherland (CS) model associated with the Hamiltonian 
\renewcommand{\theequation}{5.{\arabic{equation}}}
\setcounter{equation}{0}
\beq
H = -\frac{1}{2} \sum_{k=1}^{M}\,\left( \frac{\partial}{\partial x_k} 
\right)^2\,+\,\frac{\pi^2}{L^2} \sum_{k<l}\,\frac{\beta \left( \beta + P_{kl} 
\right)}{\sin^2 \frac{\pi}{L} \left( x_k - x_l \right)}\,, 
\label{fv.1}
\eeq
where $P_{kl}$ is the permutation operator which interchanges the `spins' of 
$k$-th and $l$-th particles, yield a realization of the standard Yangian  
algebra $Y(gl_N)$.  So it is natural to ask whether there exists some new types 
of spin CS models whose conserved quantitites would be similarly related to the 
present multiparameter deformation of the Yangian algebra.  Working along these 
lines it has been found very recently$^{29}$ that the Hamiltonians of such spin 
CS models can be easily generated by substituting the following `anyon-like' 
representation of the permutation operator $P_{kl}$ in the 
expression~(\ref{fv.1})~: 
\bea 
  &    & P_{kl}\,\left| \alpha_1 \cdots \alpha_k \cdots \alpha_l \cdots 
\alpha_M \right\rangle \nn \\ 
  &    & \qquad = \left\{ \phi_{\alpha_k 
\alpha_l}\,\prod_{\tau=1}^{N}\,\left( 
\frac{\phi_{\tau \alpha_l}}{\phi_{\tau \alpha_k}} \right)^{n_\tau} \right\}\,
\left| \alpha_1 \cdots \alpha_k \cdots \alpha_l \cdots \alpha_M \right\rangle\,,
\label{fv.2}
\eea 
where $\alpha_i\,\in\,[1,N]$ denotes the spin degrees of freedom of the $i$-th 
particle, $n_\tau$ is the number of occurrences of the particular spin 
orientation $\alpha_i\,=\,\tau$ in the configuration 
$\left| \alpha_1 \cdots \alpha_k \cdots \alpha_l \cdots \alpha_M \right\rangle$ 
when the index $i$ runs from $k+1$ to $l-1,$ $\phi_{\alpha \alpha}\,=\,\pm 1$ and 
$\phi_{\alpha \beta}$ (for $\alpha\,\neq\,\beta$) are the multiple deformation 
parameters, the above defined operator $P_{kl}$ not only interchanges the spins 
of $k$-th and $l$-th particles but also picks up a phase factor depending on 
the orientations of all the intermediate spins.  Consequently such 
representations of permutation operators would lead to novel variants of the CS 
model which contain highly nonlocal spin dependent interactions and can be 
solved exactly$^{29}$ by the action of a generalized antisysmmetric projection 
operator on the eigenvectors of the Dunkl operators.  Since a `frozen' limit of 
the spin CS model~(\ref{fv.1}) corresponds to the Haldane-Shastry (HS) spin 
chain, it is to be expected that the HS model would also admit similar 
generalization respecting the multiparametric Yangian symmetry.  Moreover, it 
may be noted that the Hamiltonian of the one-dimensional Hubbard model on an 
infinite chain respects$^{30}$ the Yangian symmetry 
$Y\left(sl_2\right) \oplus Y\left(sl_2\right)$ as an extension of the known 
$sl_2 \oplus sl_2$ slymmetry.  Thus, it will be of interest to search for a 
generalized form of the Hubbard model which would possess the multiparametric 
and coloured extension of the Yangian symmetry. 

\vspace{1.5cm}

 
\noindent{\bf References}

\medskip

\begin{enumerate}

\item R. J. Baxter, {\it Exactly Solved Models in Statistical Mechanics} 
(Academic Press, London, 1982).

\item L. D. Faddeev, {\it Int. J. Mod. Phys.} {\bf A10}, 1845 (1995).

\item P. Kulish and E. K. Sklyanin, in {\it Lecture notes in Physics} 
Vol. 151 eds. J. Hietarinta {\it et al.} ( Springer, Berlin, 1982). 

\item V. E. Korepin, N. M. Bogoliubov, and A. G. Izergin, {\it Quantum Inverse 
Scattering Method and Correlation Functions } (Cambridge Univ. Press, 
Cambridge, 1993).

\item V. G. Drinfeld, in {\it Proceedings of the International Congress of 
Mathematicians}, Berkeley, 1986 (The American Mathematical Society, Providence, 
RI, 1987) Vol. 1.

\item A. N. Kirillov, and N. Yu. Reshitikin, in {\it New Developments in the 
Theory of Knots}, Advanced Series in Math. Physics. Vol. 11, ed. T. Kohno 
(World Scientific, Singapore, 1989). 

\item M. Wadati, T. Deguchi, and Y. Akutsu, {\it Phys. Rep.} {\bf 180}, 247 
(1989).

\item V. Chari and A. Pressley, {\it A  Guide to Quantum Groups} (Cambridge 
Univ. Press, Cambridge, 1994) and references therein.

\item N. Yu. Reshetikhin, {\it Lett. Math. Phys.} {\bf 20}, 331 (1990). 

\item R. Chakrabarti and R. Jagannathan, {\it J. Phys.} {\bf A27}, 2023 (1994).

\item A. Kundu and P. Truini, {\it J. Phys.} {\bf A28}, 4089 (1995).

\item D. Bernard, M. Gaudin, F. D. M. Haldane, and V. Pasquier, {\it J. Phys.} 
{\bf A26}, 5219 (1993). 

\item K. Hikami and M. Wadati, {\it J. Phys. Soc. Japan} {\bf 62}, 4203 (1993). 

\item K. Hikami, {\it J. Phys.} {\bf A28}, L131 (1995).

\item F. D. M. Haldane, in {\it  Proceedings of the 16th Taniguchi Symposium},  
Kashikijima, Japan, 1993, eds. A. Okiji, {\it et al.} (Springer, Berlin, 
1994). 

\item A. Kundu and B. Basu-Mallick, {\it J. Phys.} {\bf A25}, 6307 (1992). 

\item L. D. Faddeev, N. Yu. Reshetikhin, and L. A. Takhtajan, in {\it 
Yang-Baxter Equation in Integrable systems}, Advanced series in Math. Phys. 
Vol. 10, ed. M. Jimbo (World Scientific, Singapore, 1990).  

\item A. Schirrmacher, {\it Z. Phys.} {\bf C50}, 321 (1991).

\item E. G. Floratos, {\it Phys. Lett.} {\bf B252}, 97 (1990). 

\item J. Weyers, {\it Phys. Lett.} {\bf B~240}, 396 (1990). 

\item R. Chakrabarti and R. Jagannathan, {\it J. Phys.} {\bf A24}, 1709 (1991).

\item R. Chakrabarti and R. Jagannathan, {\it J. Phys.} {\bf A24}, 5683 (1991). 

\item V. Karimipour, {\it J. Phys.} {\bf A26}, 6277 (1993).

\item B. Basu-Mallick, {\it Int. J. Mod. Phys.} {\bf A11}, 715 (1996).  

\item  B. Basu-Mallick and P. Ramadevi, {\it Phys. Lett.} {\bf A211}, 339 
(1996). 

\item A. N. Kirillov and N. Yu. Reshetikhin, {\it Lett. Math. Phys.} {\bf 12}, 
199 (1986). 

\item V. Chari and A. Pressley, {\it L'Enseignement Math.} {\bf 36}, 267  (1990). 
 
\item V. Chari and A. Pressley, {\it J. Reine Angew. Math.} {\bf 417}, 87 (1991).

\item B. Basu-Mallick, `Spin dependent extension of Calogero-Sutherland model 
through anyon-like representations of permutation operators', 
Preprint TIFR/TH/96-08, hep-th/9602107 

\item D. B. Uglov and V. E. Korepin, {\it Phys. Lett.} {\bf A190}, 238 (1994).  

\end{enumerate} 

\end{document}